\documentstyle[12pt]{article}
\newcommand{\vs}{\vspace{.1in}}

\newcommand{\ds}{\displaystyle}

\newcommand{\noin}{\noindent}

\def\C{I\!\!\!\!C}
\def\qed{QED. \vs}

\def\w{{\bar w}}

\def\d{{\kappa }}

\def\x{{\bar x}}

\def\k{{\bar k}}

\begin{document}

\begin{center} {\bf Uniform Zariski's Theorem On 
Fundamental Groups} \end{center}

\begin{center} {\bf Shulim Kaliman \\

Department of Mathematics and Computer Science \\

University of Miami \\

Coral Gables, FL  33124 } \\[8ex] \end{center}

\noin {\bf I. Introduction.} \vs

In [1] Zariski proved the following remarkable theorem. Let $\tilde{H}$ 
be an algebraic hypersurface in $\C I\!\! P^{n}$ where $n\geq 3$. 
Then for a generic projective plane $\tilde{A} \hookrightarrow \C 
I\!\! P^{n}$ the embedding $\tilde{A} -\tilde{H} \hookrightarrow \C 
I\!\! P^{n} -\tilde{H}$ generates an isomorphism $\pi _{1} (\tilde{A} 
-\tilde{H} )\rightarrow \pi _{1} (\C I\!\! P^{n} -\tilde{H} )$ of 
the fundamental groups. This implies the similar fact for an algebraic 
hypersurface $H$ in the Euclidean space $\C^{n}$. Consider a family 
of hypersurfaces $\tilde{H}_{p} \subset \C I\!\! P^{n}$ (resp. $H_{p} 
\subset \C^{n}$) depending algebraically on parameter $p$
from an algebraic variety $P$. It is 
natural to ask whether there exists a projective plane $\tilde{A}$ (resp. 
an affine plane $A$) such that the embedding 
$\tilde{A} -\tilde{H}_{p} \hookrightarrow \C I\!\! P^{n} -\tilde{H}_{p}$ 
(resp. $A-H_{p} \hookrightarrow \C^{n} -H_{p} $) generates an 
isomorphism of the corresponding fundamental groups for every $p$. 
In both cases the answer 
is negative, since it is enough to consider the set of all hyperplanes as this 
family of hypersurfaces depending algebraically on parameter. But the 
affine case has some advantage. Namely, we can change the coordinate system 
in $\C^{n}$ using a polynomial coordinate substitution (i.e., we can change the set of planes), whenever in 
the projective case we can use linear substitutions only. This observation 
leads to the main result of this paper. 

{\sl For every family of algebraic 
hypersurfaces $H_{p}$ in $\C^{n}$ depending algebraically 
on parameter $p \in P$,  one 
can choose a coordinate system in $\C^{n}$ in such a way that for 
some plane $A_0 \subset \C^{n}$ the embedding 
$A_0 -H_{p} \hookrightarrow \C^{n} -H_{p}$ 
generates an isomorphism $\pi _{1} (A-H_{p} )
\rightarrow \pi _{1} (\C^{n} -H_{p} )$ for each $p$. }

The scheme of the proof can be described as follows. The matter
can be reduced to the three-dimensional case.
Let $\rho : \C^3 \to \C^2$
be a projection to an $(x,y)$-plane,
let $A_0$ be $\rho^{-1} (L^0 )$ where $L^0$ is
the line $\{ y= 0 \} $ in the $(x,y)$-plane.
In section 3 we establish when the embedding $A_0 - H \subset 
\C^3 - H$ generates an isomorphism of the fundamental groups.
It turns out that the following conditions are sufficient.\vs

\noin (1) 
The embedding
$\rho ^{-1} (o) - H \subset A_0 - H$ generates an
epimorphism of the fundamental groups
where $o$ is the origin of the $(x,y)$-plane. \vs

\noin (2) 
The mapping
$\rho \mid_H : H \rightarrow \C^2$ is finite
and, if it is $l$-sheeted, then
the mapping  $\rho \mid _{H\cap A} : H\cap A \rightarrow L^{0}$ 
is also $l$-sheeted. \vs

\noin (3) There exists a nonzero polynomial
$g \in \C[x]$ such that $g(0)=0$ and for generic
$c \in \C$ the number of points where
the curve $L^c = \{ y=cg(x) \}$ meets the image
$\Gamma$ of the ramification set of the 
mapping $\rho \mid _{H} $
coincides with the number of points where
$L^0$ meets $\Gamma$. \vs

The proof of this fact is a modification of the original
argument of Zariski who dealt actually with the case
when $g(x)=x$, i.e. $\{ L^c \}$ is a pencil of lines.

The trouble with the linear coordinate substitutions is that those of them
for which one of these three conditions does not hold form a 
subvariety of
codimension 1 in the space of the linear coordinate substitutions. 
Therefore, for every family of hypersurfaces we shall construct
a wider space of polynomial
coordinate systems so that for every hypersurface 
from our family
the subset of the coordinate systems in this space
with one of these three conditions violated 
(relative to this hypersurface) has codimension at least
$l$ where the ``bigger" space of coordinate systems 
we choose the larger $l$ is. When $l$ is larger than the dimension
of the parameter set $P$ we can find a coordinate system
such that the embedding $A_0 - H_p \subset \C^3 - H_p$
generates an isomorphism of the fundamental groups for
every $p$.
 
The paper is organized as follows. After preliminaries
(section 2) we show that the embedding
$A_0 -H \hookrightarrow \C^{n} -H$ 
generates an isomorphism of the fundamental groups
provided (1)-(3) hold. In section 4
we prove some technical fact which enable us to describe
coordinate systems for which conditions (1)-(3) are true.
In the last section we prove the main result for the
three-dimensional case first and then we reduce
the general case to the three-dimensional one. \vs

\noin {\bf 2. Terminology and Notation.} \vs

In this paper $P$ and later $Q$ denote always sets of parameters, and a set 
of parameters is always an algebraic variety over the field of complex 
numbers. \vs

{\bf Definition 2.1} Suppose that $T$ is a closed algebraic subvariety in 
$\C^{n} \times P$. Let $\rho _{1} :T\rightarrow \C^{n}$ 
and $\rho _{2} :T\rightarrow P$ be the natural projections and
let $H_{p}=\rho _{1} \circ \rho ^{-1}_{2} (p)$ where $p \in P$.
Then we say that
$\{ H_p \}$ is a family of algebraic varieties with parameter
$p \in P$. 
We shall mostly deal with the case
when $H_p$ is a closed algebraic hypersurface in $\C^{n}$
for every $p \in P$. 
In this case we say that $\{ H_{p}\}$ is a
family of hypersurfaces in $\C^{n}$ with parameter $p\in P$. 
When $n=2$ we speak about a family 
of curves. \vs

%{\bf Remark 2.2} In some situations which we 
%are going to consider $\rho_2$ is,
%possibly,  not surjective. In these cases
%we shall always make it
%surjective by removing the appropriate
%subvariety of $P$. What is left is a constructive set,
%i.e. the union of algebraic varieties. By treating
%this union as a disjoint one, we
%can suppose again that the set of parameters is an
%algebraic variety.
%
{\bf Remark 2.2} We do not suppose a priori that $P$ is smooth, irreducible, 
or that $P$ has a pure dimension. But we can always 
modify $P$ to make it simpler. For instance, 
we can treat $P$ as the disjoint union
of the strata of its canonical 
stratification as a singular algebraic
variety. 
Moreover, we can modify $P$ further.  Consider, for instance,
a function $\xi$ on $P$ with a finite number of values such that
the preimage of each of these values is an algebraic subvariety of $P$.
Consider a component $P'$ of $P$ and suppose that $\xi$ is constant on 
$ P'-P'_{0}$ where $P'_{0}$ is a closed proper algebraic subvariety 
of $P'$. We can replace $P'$ by the disjoint union of $P'-P'_{0}$ and $P'_{0}$. 
Then we can replace $P'_{0}$ again by the disjoint union of the strata
of its canonical stratification. Therefore, induction by dimension implies that one can suppose from the beginning that 
$P$ is smooth and $\xi$ is constant on every component of $P$.
In particular, for every family of hypersurfaces
$\{ H_p \}$ we can suppose that the degree (or the Newton
polygon) of
the defining polynomial $f_p$ of $H_p$ is constant on every component
of $P$. Each polynomial $f_p$ can be represented
as $ \prod_{i=1}^{l(p)} (f_{p,i})^{n_i(p)}$ where
$f_{p,i}$ are irreducible polynomials which are
non-proportional for different values of $i$, and
$n_i(p)$ are natural. We can suppose again
that the function $l(p)$ is constant on every component of $P$.\vs

{\bf Convention 2.3} Since we consider the fundamental group of the
complement of a hypersurface it is natural to restrict ourselves to
the case when this hypersurface is reduced. In the case
of a family of hypersurfaces we can require that each
generic hypersurface from this family is reduced.
In general, even
if a generic member of a family of hypersurfaces
is reduced there may be some hypersurfaces which are not reduced.
But dividing the set of parameters into a disjoint union
of subsets in the manner described in
Remark 2.2, we can suppose and we {\sl will suppose}
that every hypersurface in our
family is reduced (for this we have to suppose that
$l(p)$ from Remark 2.2 is constant on every component
of $P$ and to replace
$f_p= \prod_{i=1}^{l(p)} (f_{p,i})^{n_i}$ by
$\prod_{i=1}^{l(p)} f_{p,i}$ on each component of $P$). \vs

{\bf Remark 2.4} We say that an algebraic
variety $B$ depends algebraically on a polynomial $f$
if $B$ can be viewed as an member
of a family of algebraic varieties whose parameter set
is a Zariski open subset of the set of polynomials
of fixed degree (or with a fixed Newton polygon). \vs

{\bf Definition 2.5} Let $H$ be a reduced algebraic 
hypersurface in $\C^{n}$ and $A$ be a closed 
affine algebraic subvariety in $\C^{n}$. We say that $A$ is $H$ 
compatible if the 
embedding $i:A-H\hookrightarrow \C^{n} -H$ generates an 
epimorphism of the fundamental 
groups $i_{*} :\pi _{1} (A-H)\rightarrow \pi _{1} 
(\C^{n} -H)$. We say that $A$ is strictly $H$ compatible if $i_{*}$ is 
an isomorphism.  We say that $A$ is (strictly) compatible
relative to a family of hypersurfaces if it is (strictly)
$H$ compatible for every hypersurface $H$ in this family.  \vs

%{\bf Definition 2.6} If a closed affine algebraic subvariety $A$ in 
%$\C^{n}$ is 
%isomorphic to $\C^k$
%we call it a $\C^k$-manifold. When $k=1$ we call 
%$A$ a $\C$-curve and when $k=2$ we call
%it a $\C^2$-surface. \vs
%
%Our aim is to prove that for every family of 
%hypersurfaces $\{ H_{p}\}$ in $\C^{n} \;\; (n\geq 3)$ with 
%parameter $p\in P$ there 
%exists a surface $A_{0}$ which is strictly $H_{p}$ 
%compatible for every 
%$p\in P$ and  which becomes a coordinate plane after a polynomial
%coordinate substitution. \vs

{\bf Definition 2.6} We say that some property holds for a generic point 
of an algebraic variety $P$ if for every irreducible component $P'$ of $P$ 
there exists a proper closed
algebraic subvariety $P'_{0}$ such that this property 
is true for every $p\in P'-P'_{0}$. \vs

In particular if we have several properties and each of them holds for a 
generic point of $P$, then all of them hold simultaneously for a generic 
point of $P$. \vs

%We conclude this section with citing notation we shall use in the
%remainder of this article. We denote
%by $G_{m}$ the space of 
%polynomials in one variable of degree 
%$\leq m$. \vs 
\newpage

\noin {\bf 3. Modification of 
Zariski's approach for the non-linear case.} \vs

In this section the 
projection $\rho :\C^{3} \rightarrow \C^{2}$ is
given by  $(x, y, z)\rightarrow (x, y)$ 
and $H\subset \C^{3}$ is a reduced hypersurface
whose defining polynomial is $f$. Denote by $A_0$ the $(x,z)$-plane
in $\C^3$ and by $L^0$ the $x$-axis in the $(x,y)$-plane (of course,
$A_0 = \rho^{-1}(L^0)$). Suppose that $g \in \C [x]$
is a nonzero polynomial such that $g(0)=0$. \vs

{\bf Definition 3.1} We say that the $\C$-curve $L^{0}$ 
is appropriate with respect to the triple
$(H, \rho ,g)$ if the following properties hold \vs

\noin (1)
The $z$-axis (i.e., the line
$\rho ^{-1} (o) $
where $o$ is the origin of the $(x,y)$-plane)
is $H \cap A_0$ compatible in the plane $A_0$. \vs

\noin (2) 
The mapping
$\rho \mid_H : H \rightarrow \C^2$ is finite
and, if it is $l$-sheeted, then
the mapping  $\rho \mid _{H\cap A_0} : H\cap A_0 \rightarrow L^{0}$ 
is also $l$-sheeted. \vs

\noin (3) 
For generic $c \in \C$ the number of points where
the curve $L^c := \{ y=cg(x) \}$ meets the image
$\Gamma$ of the ramification set of the 
mapping $\rho \mid _{H} $
coincides with the number of points where $\Gamma$ meets $L^0$. \vs

The aim of this section is \vs

{\bf Theorem 3.2} {\em Let $L^{0}$ be appropriate with respect to $(H, \rho
,g)$. Then $A_{0}$ is strictly $H$ compatible.} \vs

In the case when $g$ is linear
this theorem can be extracted from the original paper of Zariski [1]. In 
general case we follow also the outline of his arguments.
The proof consists of several lemmas and
we discuss first some {\sl additional notation } which will be used
further in these lemmas.

Let $X$ be the set of points
where $L^0$ meets $\{ g (x)=0 \}$.
The family $\{ L^c \}$ can be viewed as a linear system of curves 
whose base point set is $X$.
This linear system of curves generates the 
mapping $\tau :\C^{2} -X\rightarrow \C I\!\! P ^{1}$ such that 
$\tau ^{-1} (c)=L^{c}-X$ for $c\in \C I\!\! P^1$
(where $L^{\infty }=\{ g(x)=0\}$). 
There exists a finite set $C=\{ c_{1} ,..., c_{r} \} \subset \C$ 
such that for every $c\in \C -C$ the curve 
$\tau ^{-1} (c)$ meets $\Gamma $ at the same number of points.

Put $Z = X \cup \tau^{-1}(C \cup \infty )$ 
(i.e. $Z= L^{c_1} \cup \ldots \cup L^{c_r} \cup L^{\infty }$).
Set ${\Gamma}' = \Gamma -Z$ (i.e. ${\Gamma }' = 
\Gamma \cap \tau^{-1} (\C
-C)$) and
${\Gamma}'' = \rho^{-1}(
{\Gamma }') \cap H$.
The set of points
where a generic curve $L^c$ meets $\Gamma$ non-normally
is contained in $X$, by condition (3). 
Hence every non-smooth
point of $\Gamma$ is contained in $Z$ since at these points
none of generic $\C$-curves $L^c$ can meet $\Gamma$ normally.
Hence ${\Gamma }'$ is smooth and we need to
show that the mapping $\rho \mid_{{\Gamma }''} : {\Gamma }''
\to { \Gamma }'$ is unramified. This fact can checked locally. \vs

{\bf Lemma 3.3} {\sl Let ${\cal H}$ be the germ of an analytic
surface at the origin of $\, \C^3$. Suppose that the mapping
$\rho_0 : {\cal H} \to (\C^2 ,o)$ is finite
where $\rho_0 =\rho \mid_{\cal H}$.
Let $\gamma \subset (\C^2, o)$ be the image of the ramification set
for the mapping $\rho_0$ and let
${\cal K} = \rho^{-1}_0 (\gamma )$. Suppose that
$\gamma$ is smooth and $\rho_{\cal K} : {\cal K}
\to \gamma$ is the restriction of $\rho$ to ${\cal K}$.
Then $\rho_{\cal K}$ is unramified. } \vs

{\sl Proof.}
Let ${\cal L}$ be the preimage of
a generic point $b
\in (\C^2, o) - \gamma $ under the mapping $\rho_0$.
The fundamental group  of
$(\C^2 ,o) - \gamma$ is isomorphic to the group of integers
since $\gamma$ is smooth and this group
acts on ${\cal L}$. Hence ${\cal L}$
can be represented as the disjoint union of minimal invariant
subsets of ${\cal L}$
relative to this action. These subsets
correspond to the irreducible components of ${\cal H}$.
Let $s$ be the number of points
in the preimage ${\cal L}_0$ of a generic point $a \in
\gamma$ under $\rho_0$.
Sending $b$ to $a$ one can see that the points
of ${\cal L}_0$ generate a partition of ${\cal L}$ into
disjoint subsets. Using the fact that the fundamental group of
$(\C^2 ,o) - \gamma$ is the group of integers,
one can check that this partition is invariant under
the action of the fundamental group (it is enough
to consider the action on ${\cal L}$ of a small simple loop
from $(\C^2 ,o) - \gamma$ around $a$ since this loop
can be viewed as a generator of the fundamental group).
Hence we can represent
${\cal H}$ as the union $\bigcup_{i=1}^s {\cal H}_i$ of
the germs of surfaces
so that ${\cal H}_i \cap \rho^{-1} (\gamma ) = {\cal K}_i$
where ${\cal K}_i \subset {\cal K}$ are the germs of different curves
with the origin as the only common point (since
${\cal K}_i \cap \rho^{-1} (a)$ is exactly one point
in ${\cal L}_0$ that corresponds to ${\cal H}_i$). Assume that
$s \geq 2$. Consider the germ of the curve $\zeta = {\cal H}_1
\cap {\cal H}_2$.  By construction, the germ
$\rho (\zeta )$ meets $\gamma$ at the origin only.
But it must be contained in $\gamma$ since $\zeta$ is
contained in the ramification set. Contradiction.
Hence $s=1$ which is the desired conclusion. \\
\qed 

Put $Y = \tau^{-1} (\C - C) \cup X$. For every subset $K$ of $\C^2$
we denote by $\tilde {K}$ the set $\rho^{-1} (K) -H$. For instance,
$\tilde {Y} = \rho^{-1} (Y) -H$. \vs

{\bf Lemma 3.4} {\sl
The natural embedding $i:\tilde {Y} \hookrightarrow \C^{3} -H$ 
generates an isomorphism of the fundamental groups $i_{*} 
:\pi _{1} (\tilde {Y}) \rightarrow \pi _{1} (\C^{3} -H)$.} \vs

{\sl Proof.} By condition (3) in Definition 3.1, 
$C$ does not contain 0. Since 
$\tau ^{-1} (0)=L^{0} -X$ the set $\tilde {Y}$ 
contains $A_{0} -H$. We shall see later (Lemma 4.2)
that under condition (2) every line which is parallel
to the $z$-axis and which meets $H$ at $l$ points, is
$H$ compatible. Thus condition (2) implies that
$A_0$ is $H$ compatible.
% and, therefore, $\rho^{-1}(o)$ is $H$ compatible by
%condition (1). 
Hence $i_{*}$ is 
an epimorphism. 

Let $\delta _{1} ,..., \delta _{s}$ be a set of 
generators in $\pi _{1} (\tilde {Y})$, and, therefore, it can be 
treated as a set 
of generators in $\pi _{1} (\C^{3} -H)$. We need to show that if 
$\delta _{1} ,..., \delta _{s}$ satisfy some generating relations in 
$\pi _{1} (\C^{3} -H)$ then
they satisfy the same relations in $\pi _{1} (\tilde {Y})$. For 
this purpose it suffices to show that for every 2-cell 
$\Delta _{1}$ in $\C^{3} -H$ with a boundary 
$\partial \Delta _{1} \subset \tilde {Y}$ 
there exists a 2-cell $\Delta _{2} \subset \tilde {Y}$ with 
the same boundary.
%Note that $\tilde {Y} \cap \tilde {Z}$ 
%coincides with  $\tilde{X}$. 
We can suppose 
that if $u$ is an intersection point of $\tilde {Z}$ 
and the interior of $\Delta _{1}$ then $\Delta _{1}$ meets $\tilde {Z}$ 
normally at $u$, and $u \notin \tilde {X}$. (We can do this since the
real codimension of $\tilde {X}$ in $\C^3 - H$ is 4.)
Since $X$ is the base point set for $\{ L^c \}$
one can choose a path $\xi$ in $\tilde {Z}$ 
joining $u$ and a point $v\in \tilde{X}$ in such 
manner that $\rho (v)$ is the only point from $\rho (\xi )$ that belongs to 
$X$ and that $\xi - v$ is contained in
the smooth part $\tilde {Z}^*$ of $\tilde {Z}$. One can
identify a neighborhood of $\tilde {Z}^*$ in $\C^3$ with
a neighborhood of the zero section of the normal bundle to
$\tilde {Z}^*$. Moreover, we can suppose that the intersection of
$\Delta_1 $ with this neighborhood is contained in a fiber of this bundle.
Choose a small 2-cell $\Delta _{\varepsilon } (u)\subset \Delta _{1}$ 
with center at $u$ and replace it by a cone in 
$\tilde {Y}$ with the following properties: 
$v$ is the vertex of the cone and the only point where 
the cone meets $\tilde {Z}$, 
the base of the cone coincides with the boundary of 
$\Delta _{\varepsilon } (u)$, the intersection of the cone
with the fiber of the normal bundle to $\tilde {Z}^*$
at every point of $\xi -v$ is a circle. Repeating 
this procedure we obtain a 2-cell ${\Delta }_{2} \subset \C^{3} -H$ 
such that its boundary coincides with the 
boundary of $\Delta _{1}$ and its interior meets $\tilde {Z}$ 
only at points from $\tilde{X}$. In 
particular, ${\Delta }_{2} \subset \tilde{Y}$. \\
\qed

Let $\sigma_1, \ldots , \sigma_r$ be a bouquet of simple
loops in $\C - C$ with one common point at the origin
so that these loops generate the fundamental group of
$\C -C$. Put $B = \tau^{-1} (\bigcup_k \sigma_k ) \cup X$.
We have to show that $\tilde {B}$ is a deformation retract
of $\tilde {Y}$, and, in particular, the embedding $j: \tilde {B}
\hookrightarrow \tilde {Y}$ generates an isomorphism 
$j_{*} :\pi _{1} (\tilde {B} )\rightarrow \pi _{1} (\tilde {Y})$. 
The construction of this
deformation can be reduced to a simpler problem due to \vs

{\bf Lemma 3.5} {\sl Let 
${K}$ be a subset of ${Y}$. 
Suppose that $\d = \{ \d_t \mid t \in [0,1] \} $ is a 
path in the space of continuous mappings
from $K$ to $Y$ such that 
$\d_0$ is the identical embedding,
the restriction of $\d_t$ to $K 
\cap X$ is the identical embedding
for every $t \in [0,1]$,
$\d_t ({\Gamma } \cap {K}) \subset \Gamma$, and
$\d_t (K - (X \cup \Gamma ))
\subset Y- (X \cup \Gamma )$.
Suppose also that the restriction of $\d$ to
$({K} -X) \times [0,1]$ is smooth.
Then there exists a deformation $D = \{ D_t \mid t \in [0,1] \}$ 
of the identical embedding $D_{0}: 
\tilde {K} \hookrightarrow \tilde{Y}$ such that $\rho D=\kappa \rho$ 
and $D$ is identical on $\tilde {K} \cap \tilde{X}$. Moreover,
if for every $t$ the mapping $\d_t$ is a homeomorphism
between $K$ and $\d_t (K)$ then $D_t$ is a homeomorphism
between $\tilde {K}$ and $D_t (\tilde {K})$.} \vs

{\em Proof.}
Fix a neighborhood 
$U$ of $H- \rho^{-1} (X)$ in $\C^3$
so that $ \bar {U} \cap \rho^{-1}(b)$ is compact
for every $b \in \C^2$ 
($\bar {U}$ is, of course, the closure of $U$ in $\C^3$).
Suppose also that 
$\bar {U} \cap \rho^{-1} (X) = H \cap \rho^{-1} (X)$, i.e this set
is finite.
Consider the natural projection
$T\C^2 \to \C^2$ where
$T\C^2$ is the tangent bundle of $\C^2$, and
the mapping $\rho : \C^3 \to \C^2$. 
They generate the set 
$T = \C^3 \otimes_{\C^2} T\C^2$ with  the natural projections
$pr_1 : T \to \C^3$ and $pr_2 : T \to T\C^2$. 
Put $Y'= \C^3 - (H \cap \rho^{-1} (X \cup \Gamma))$
and put $W = pr_1^{-1} (Y')$. Let
$W_0 = pr_1^{-1} (\tilde {\Gamma }') \cap pr_2^{-1} (T{\Gamma }')$ and
$\overline {W}_0 = pr_1^{-1} (\rho^{-1} ({\Gamma }'))
\cap pr_2^{-1} (T{\Gamma }')$ where $T{\Gamma }'$ is the
tangent bundle to ${\Gamma }'$.
Using partition of unity one can construct a smooth mapping
$\chi \colon W \to T {Y}'$ where $T {Y}'$ is the tangent 
bundle of ${Y}'$ with the following 
properties:

(i) for every $w \in W $ we have $ \rho_* \chi = pr_2$;

(ii) for every $w \in W $ with $pr_1 (w) \notin {U}$ the 
$z$-coordinate of the vector $\chi (w)$ is zero;

(iii) for every $w \in W$ with  $pr_1(w) \in H$ 
the  vector $\chi (w)$
is tangent to $H$. (Note
that this tangent vector exists since the 
restriction of $\rho$ to $H- \rho^{-1} ({\Gamma })$ is
an unramified covering of $\C^2 - \Gamma$);

(iv) the restriction of
$\chi $ to $W_0$ can be extended to a smooth
mapping $\chi_0 : \overline {W}_0 \to T\C^3$ 
(where $T\C^3$ is the tangent bundle of $\C^3$) so that
for every $w \in \overline {W}_0$
with $pr_1 (w) \in {\Gamma }''$ the vector
$\chi_0 (w)$ is tangent to ${\Gamma }''$.
(This vector $\chi_0 (w)$ exists
since the restriction of $\rho$ to ${\Gamma }'' $ is an 
unramified covering of ${\Gamma }'$, by Lemma 3.3.) \vs

It is worth mentioning that we need 
condition (iv) separately from (iii) since  in general the mapping $\chi$
cannot be extended continuously to the 
points $w$ with $pr_1 (w) \in \rho^{-1} ( \Gamma ) \cap H$. 

Consider the curve $ \d (b) =\{ \d_t(b) \mid t \in [0,1] \}$ for 
$b \in K-X$.
Suppose for simplicity that
it has no selfintersection points (otherwise
we can replace the curve with its graph). Then $\widetilde {\d (b)}$ 
is a smooth  real manifold.
For each $u_0 \in \widetilde { \d (b)}$ 
there exists $t_0 \in [0,1]$ such that $\rho (u_0) = 
\d_{t_0} (b)$. Let
$v_b (t_0)$ be the vector tangent to $\d (b)$ at  $\d_{t_0}(b)$ 
which is generated
by differentiation with respect to $t$. Then the vector
$\chi (u_0 \otimes v_b (t_0))$ is tangent to $\widetilde { \d (b) }$ at $u_0$. 
Therefore, such vectors 
define a vector field
on $\widetilde {\d (b)}$ . Let $\rho (u) =b$ and
let $D(u) = \{D_t (u) \mid t \in [0,1] \}$ 
be the integral curve of this vector field such
that it begins at the point $u=D_0 (u)$. These curves 
define a deformation D
of $\tilde {K} - \tilde {X}$ in $\tilde {Y} - \tilde {X}$
with $\rho D = \d \rho $ 
unless for some $u \in \tilde {K} - \tilde {X}$ the
curve $D (u)$ goes either to infinity or to $H$ for a finite time. 
It cannot go infinity for a finite time
within $U$ due to the fact that $\bar {U} \cap \rho^{-1} (b)$ is compact 
for every $b$ and the mapping
$\rho \mid_H$ is finite. 
Outside $U$ it cannot go to infinity as well, by (ii). 
When $b \notin \Gamma$ the curve $D(u)$ cannot reach $H$ for a
finite time due to (iii). When $b \in \Gamma$ the curve
$\d (b) \subset \Gamma$, by the assumption of Lemma,
and $D(u)$ cannot reach $H$ again
for a finite time due to (iv). 
Therefore, $D(u) \subset \tilde {Y}$ and
$D_t(\tilde {K} - \tilde {X}) \subset \tilde {Y}$
for every $t$.

Note that the combination of (ii) and the facts that 
$\d_t \mid _{{K} \cap {X}}$ is the identical embedding and
that the set $\bar {U} \cap \rho^{-1} (X)$ is finite implies that $D$ can 
be extended to 
$\tilde {X}$ by the identical deformation which proves the
first statement of Lemma. The second statement follows obviously from
the construction of $D$.\\
\qed

The next lemma is almost the exact repetition of Lemma 3.5
but we give its proof for the sake of completeness. \vs

{\bf Lemma 3.6} {\sl Let ${K} \subset
{Y}$ and $\d_{0} : {K} \rightarrow Y$ be the identical embedding. 
Let $\theta_0$ be the identical embedding of 
$\tau (K - X)$ into $\C -C$. Suppose that 
there exists a smooth 
path $\theta =\{ \theta _{t} | t \in [0,1] \}$ in the space
of smooth mappings from $\tau (K - X)$ to $\C -C$. 
Then there exist a deformation $\d=\{ \d_{t} \}$ of $\d_{0}$ so that
$\d_t \mid_{{X}}$ is the identical embedding for every $t$, 
$\d ({\Gamma } \cap {K} ) \subset \Gamma$, $\d 
({K} - (X \cup \Gamma ) )
\subset Y - (X \cup \Gamma )$, and $\tau \circ \d |_{K-X}
= \theta \circ \tau |_{K-X}$. Moreover, if for every $t$ the mapping
$\theta_t$ is a diffeomorphism between $\tau (K - X)$ and its image
then $\d_t$ is a homeomorphism between $K$ and its image, and
the restriction of $\d$ to $(K-X) \times [0,1]$ is smooth. } 
\vs

{\em Proof.} 
Note that $Y-X = 
\tau^{-1} (\C-C)$
can be treated as $(\C - \{ x_1^0, \ldots , x_{l}^0 \} ) 
\times (\C - C) \subset \C \times (\C - C)$
where $x_1^0, \ldots , x_{l}^0$ are the $x$-coordinates of the points
of $X$.
Then $\tau $ can be treated as the natural 
projection $\tau_2$ to the second factor. Let
$\tau_1 \, \colon \, \tau^{-1} (\C - C) \to \C$ 
be the projection to the first factor. Choose a
small tubular neighborhood $U$ of 
${\Gamma }'$ in $\tau^{-1} (\C- C)$ such that
its closure $\bar {U}$ in $\C \times (\C - C)$ 
does not meet the sets $x_i^0 \times (\C -C), \,
i=1, \ldots , l$.
Using partition of unity we can construct a vector field 
$\mu$ on $\C \times (\C -C)$ so that

- outside $U$ we have $\tau_{1*} (\mu ) \equiv 0$,

- $\mu$ is tangent to ${\Gamma }'$, 

- $\tau_{2*} (\mu )$ is a nonzero constant vector field 
on $\C - C$ (this means that
one can suppose that the phase flow associated 
with $\mu $ transforms $L^{c}$ 
into $L^{c+t}$ for time $t$  whenever this flow  is 
defined correctly). \vs

If $b\in  K \cap {X}$ put $\d (b)=b$ and 
if $b\in {K} - {X}$ define $\d_{t} 
(b)$ as follows. 
Consider $c = \tau (b)$ and 
$M_c = \tau^{-1} (\theta (c))$. Suppose for simplicity 
that $\theta (c)$
has no selfintersection points. Then $M_c$ is a smooth 
real manifold which is
naturally embedded in $\tau^{-1} (\C -C)$. For each 
$a \in M_c$ there exists $t$ such that $\tau (a) = \theta_t (c)$. 
Consider the vector $\theta '_{t} (c) \mu (a) $ at $a$ 
where ${\theta }'_t (c)$ is 
the derivative of  the function $\theta (c) \, : \, [0,1] \to \C$ 
with respect to $t$. 
This vector is tangent to $M_c$ and,
therefore, such vectors define a vector field $\mu_c$ on $M_c$.
This vector field  
defines uniquely an integral curve $\d_{t} (b)$ in $M_c$
which begins at $b = \d_0 (b)$.
 
The continuity of $\d$ is clear unless for some $b \in K -X$ the 
curve $\d (b)$ goes either to infinity or to
$x_i^0 \times (\C - C)$ for a finite time. 
But it cannot go to infinity
since we cannot reach infinity
within $U$ for a finite 
time  due to the description of $U$,
and outside $U$ the 
behavior of $\d_{t} (b)$ is defined by the vector field $\mu$ 
which does not send points from $\tau^{-1} (\C - C)$ to infinity
since $\tau_{1*} (\mu ) = 0$.
Similarly, a set
$x_i^0 \times (\C -C)$) cannot be reached for a finite time. Thus $\d$ 
is continuous and $\d_t ({K} - {X})\subset Y-X$ for every $t$. 
Note that $\mu_c$ is tangent to $\Gamma \cap M_c $. Hence the curve 
$\d (b)$ is either contained in
$\Gamma $ or does not meet it which yields the desired properties of $\d$
in the first statement. The second statement follows obviously from the
construction of $\d$ and the fact that the integral curve
$\d (b)$ depends smoothly on $b \in K - X$.\\ 
\qed

{\sl The Proof of Theorem 3.2.} We can suppose that the
bouquet $\bigcup_k \sigma_k$ is chosen so that there exists
a smooth deformation $\theta$ of $\C -C$ to this bouquet.
By Lemma 3.6, there exists a 
deformation
$\d =\{ \d_{t} \}$ of $Y$ to $B$ with the
prescribed properties. 
By Lemma 3.5, there exists a 
deformation of $\tilde {Y}$ to $\tilde {B}$
which is identical on $\tilde {X}$.
Thus it remains to show that the natural 
embedding of $A_0 -H$ into $\tilde {B}$ generates
an isomorphism of the fundamental groups. 
As we mentioned in the proof of Lemma 3.4 this
embedding generates an epimorphism of the fundamental groups.
By (1) the line $\rho^{-1} (o)$ is $H$ compatible.
Hence the loops $\delta_1 , \ldots , \delta_s$ which generate
$\pi_1 (\rho^{-1} (o) -H)$ can be viewed as generators of
$\pi_1 (\tilde {B} )$. If we travel along a loop $ \sigma_j
: [0,1] \to \C - C$
from the point $o= \sigma_j (0)$ to a point $c =\sigma_j (t)$ 
then Lemma 3.6 provides us
with an appropriate homeomorphism between $L^0$ and $L^c$
which depends continuously on $t$.
It generates in turn a homeomorphism
(depending continuously on $t$) between
$\tilde {L}^0 \, (=A_0 -H)$ and $\tilde {L}^c$ which is
identical on $\tilde {X} = \tilde {L}^0 \cap \tilde {L}^c$,
by Lemma 3.5. Thus after traveling along the whole loop
$\sigma_j$ we deform each element $\delta \in \pi_1 (\tilde {L}^0)$
into another element $\delta^{(j)} \in \pi_1(\tilde {L}^0)$.
By the van Kampen theorem [3], the generating relations
for $\delta_1, \ldots , \delta_s$ in $\pi_1 (\tilde {L}^0)$
together with the relations $\delta_i = \delta^{(j)}_i$
give all the generating relations between $\delta_1 , \ldots ,
\delta_s$ in $\pi_1 (\tilde {B})$. But since
for every $j$  the homeomorphism
of $\tilde {L}^0$ generated  by the loop $\sigma_j$ is identical
on $\tilde {X}$ the loops $\delta_i$ and $\delta^{(j)}_i$
coincide for every $i$. Hence the
embedding of $\rho^{-1} (o) -H$ into $\tilde {B}$ generates
an isomorphism of the fundamental groups. Since
$\tilde {B}$ is a deformation retract of $\tilde {Y}$,
Lemma 3.4 implies that the
embedding of $A_0 -H$ into $\C^3 -H$ generates also
an isomorphism of the fundamental groups. \\
\qed

\noin {\bf 4. Technical Facts } \vs

In this section  we shall describe polynomial coordinate
substitutions which lead eventually to conditions (1)-(3)
from Definition 3.1.
The fist two lemmas of this section are formulated
for hypersurfaces in $\C^3$ but they can be easily reformulated
for hypersurfaces in $\C^n$ with $n \geq 3$. \vs

{\bf Lemma 4.1} {\em Let $(x_1, y_1 , z_1)$ be a coordinate 
system in $\C^{3}$ and let
$\{ H_p \subset \C^{3} \mid p \in P \}$ be a family of 
hypersurfaces with defining polynomials $\{ f_p \}$. Suppose that
$d$ is natural such that $d > {\ds \max_{p \in P}}
\, deg f_p$. 
Let $(x,y,z)$ be a new coordinate system such that
$(x_1,y_1,z_1) = (g_1(x,y,z), g_2(x,y,z),z)$ where
$g_1$ and $g_2$ are polynomials.
Suppose that the degrees of
$g_1(a,b,z)$ and $ g_2(a,b,z)$ are
$d_1$ and $d_2$ respectively for every constants $a$ and $b$.
Let $d_2>d$ and $d_1 >d \cdot d_2$. 
Then the restriction of the projection $\rho (x,y,z) \to (x,y)$
to every hypersurface $H_p$
is a finite morphism.}

{\em Proof}.  Fix $p \in P$.
The restriction of $\rho$ to $H_p$ is finite if
for every constants $a$ and $b$
the degree of the polynomial 
$\varphi (z) = {f}_p (g_1(z,a,b), g_2(z,a,b),z)$ 
does not depend on $a$ and $b$.
Consider monomials
$x_1^{k} y_1^{l} z_1^m$ which
are present in $f_p$ with nonzero coefficients and 
consider the vectors $(k,l,m)$.
Suppose that $(k^0, l^0 , m^0)$ is the 
greatest among these vectors in the lexicographic order.
One can see that $deg \, \varphi$ coincides then 
with $k^0 d_1 + l^0d_2+ m$
regardless of the choice of $a$ and $b$.\\
\qed

{\bf Lemma 4.2} {\em Let
$\rho : \C^{3} \to \C^{2}$ be the projection given by
$(x,y,z) \to (x,y)$.
Suppose that a reduced hypersurface $H$  in $\C^3$
does not contain lines parallel to the $z$-axis. 
Consider a line $L=\rho^{-1} (w)$ with $w \in \C^{2}$ which 
meets $H$ at $deg_{z} f$ points (counting without multiplicity). Then 
$L$ is $H $ compatible.} \vs

{\em Proof.} There is an algebraic subvariety
$S\subset \C^{2}$ such that the line $\rho^{-1}(s)$ 
meets $H$ at less than $deg_{z} f$ points counting 
without multiplicity iff $s \in S$. 
Put $E=\C^{3} - (H \cup \rho ^{-1} (S))$ and 
$\tau =\rho \mid _{E}$. Then the mapping 
$\tau :E\rightarrow \C^{2} -S$ is a fibration whose generic 
fiber $F$ is a $deg _{z} f$ times 
punctured complex line. Let $i:F\rightarrow E$ be the natural embedding. 
We have the exact sequence of the fundamental groups 
$$\rightarrow 
\pi _{1} (F) \stackrel{i_{*}}{\rightarrow } \pi _{1} (E) 
\stackrel{\tau _{*}}{\rightarrow } \pi _{1} (\C^{2}
-S)\rightarrow 0.$$ 
Choose simple loops $\{ \sigma _i \}$ in $\C^{2} - S$ around
each point in a finite set $\{ s_i \} \subset S$  such that these
loops generate the whole fundamental group
$\pi_1 (\C^{2}  - S)$. Since $\rho^{-1}(s_i)$ is 
not contained in $H$ one can choose
a loop $\gamma_i$ in $E$ so that it is contractible in 
$\C^{3} - H$ and $\tau (\gamma_i ) = \sigma_i$.
Then the exact sequence 
implies that every element in $\pi _{1} (E)$ 
can be written in the form $uv$ where 
$v\in i_{*} (\pi _{1} (F))$ and $u$ lies in 
the group generated by $\{ [\gamma _{i} ] \}$. 
Consider the embedding $j:E\rightarrow 
\C^{3} -H$ and the corresponding homomorphism of the fundamental 
groups $j_{*} :\pi _{1} (E) \rightarrow \pi _{1} (\C^{2} -H )$.
Note that $[\gamma _{i} ]\in ker \, j_{*}$ and $j_{*} (uv)=j_{*} (v)$. 
On the other hand $j_{*}$ is surjective, of course. 
Hence $j_{*} \circ i_{*} (\pi _{1} 
(F))=\pi _{1} (\C^{3} -H )$ and we are done, 
since $L-H $ is a generic fiber of $\tau$.\\
\qed

We shall cite notation which will be used in the
remainder of this section. A curve $\Gamma$ in $\C^{2}$ is 
always reduced and 
it coincides with the zero locus of a 
polynomial $f$.
We shall denote by $h^0$ a polynomial
in one variable of degree $d_0$.
We shall consider $\C$-curves in $\C^{2}$ given by
equations of form $y+h^{0} (x)+g(x)=0$ 
where $g$ runs over the space of polynomials $G_{m}$ 
of degree $m<d_0$. 
The family of these curves will be denoted by $V_{m} 
(h^{0} )$.
There is a natural bijection between $G_{m}$ and 
$V_{m} (h^0)$ and we denote by $L_{g}$ 
the $\C$-curve from 
$V_{m} (h^{0} )$ that corresponds to $g\in G_{m}$. \vs

{\bf Lemma 4.3} 
{\sl Let $\Gamma$ be an algebraic curve in $\C^2$ which does 
not contain lines parallel to the $y$-axis.
Let a polynomial $g$ has a simple root $x^{0}$ so that 
the line $C = \{ x=x^0 \}$
meets $\Gamma$ at $deg _{y} f$ different points.
Let $g^0 \in G_{m}$ and $g(c) = g^0+ cg$.
Then the curve $L_{g(c)} \in V_{m}
(h^0)$ is $\Gamma$ 
compatible when $|c|$ is sufficiently large.} \vs

{\em Proof.} By Lemma 4.2, $C$ is $\Gamma$ compatible. 
When $|c| \rightarrow \infty$ then $L_{g(c)}$ approaches to $C$.
Repeating the argument of [2, Lemma 3], one can see
that $L_{g(c)}$ is $\Gamma$
compatible for large $\mid c\mid$. (In [2, Lemma 3] we used
the term ``$\Gamma$ proper" instead of ``$\Gamma$ compatible". We made
this replacement since the term ``proper" may be
misunderstood.) \\
\qed

{\bf Remark 4.4} Suppose that $\omega \in G_{m}$
is a nonzero polynomial and $g^0 \in G_{m}$.
Let $V$ be an affine subspace of $V_{m} (h^0)$
so that it consists of the $\C$-curves of
form $\{ y+h^0(x)+ g^0 (x)+
\omega (x) \tilde {g} (x)=0 \}$ where 
$\tilde {g} \in G_{n}$ and 
$ n=m-deg \omega$.
Note that if $n >0$ then
for generic $\tilde {g} \in G_{n}$
the polynomial $g= \tilde {g} \omega$ has always a simple 
root  $x^0$ such that the 
line $C$ described in Lemma 4.3 is $\Gamma$ compatible. \vs

{\bf Lemma 4.5} {\em Suppose that
the curve $L_{{g}^{0}}$ 
meets $\Gamma $ at the same number of
points as a generic 
$\C$-curve in $V$ where $L_{g^0}$ and
$V$ is from Remark 4.4 with $n>0$. 
Suppose also that
$\Gamma$ does not contain lines parallel to the $y$-axis.
Then $L_{{g}^{0}}$ 
is $\Gamma$ compatible.} \vs

{\em Proof.} Put ${\cal L} = \{ ((x,y) ,L) \in \C^{2} \times V
\mid (x,y) \in L -\Gamma \}$. Let $\kappa : {\cal L} \to V$ be
the natural projection. Then there exists a closed algebraic
subvariety $S$ of $V$ such that
a curve $L$ from $V$ meets $\Gamma$ at
a less number of points than the generic curve from $V$
iff this curve $L$ is from $S$. Hence the restriction of
$\kappa$ to $\kappa^{-1}(
V -S)$
is a fibration over  
$V -S$.
This fibration provides an isotopy between 
$L_{g^0} -\Gamma$ and 
$L-\Gamma $ in $\C^{2} -\Gamma$ where $L$ is
a generic curve in $V$.
By Lemma 4.3 and Remark 4.4, $L$ and, therefore,
$L_{g^0}$ are $\Gamma$ compatible. \\
\qed

The number of points at which $L_g$   meets
$\Gamma$ may change when $L_g$
runs over $V_{m} (h^{0} )$ 
but at least we can fix the number of intersection points of 
$L_g$ and $\Gamma$ counting with multiplicity. \vs

{\bf Lemma 4.6} {\sl  Let $\{ \Gamma_p \}$ be a family
of curves in $\C^{2}$ with parameter $p \in P$ and let
$d_0$ be natural so that $d_0 > m +k$ where $k= dim \, P$ and $m>0$.
Then for generic $h^0 \in G_{d_0}$
and every $p \in P$ the intersection number
$L_g \cdot \Gamma_p$ is finite and it
does not depend on $g \in G_{m}$. In particular,
$L_g$ is not a component of $\Gamma_p$ for every $p$ and $g$.} \vs

{\sl Proof.}
Let $N$ be the maximal possible degree of the 
polynomial $ \varphi (x) =f(x, {h}(x))$
where ${ h}$ runs over $G_{d_0}$. 
Let $L$ be the curve $y = h(x)$.
Suppose that $L \cdot \Gamma = N-l$. Then
the degree of $f(x,h(x)) = \sum^N_{i=0} a_ix^i$ is $N-l$.
Note that the coefficients  $a_i=0$ for $i > N- l$ only when
$h$ belongs to an algebraic subvariety 
${\cal A}(l)$ of $G_{d_0}$ which depends
algebraically on $f$. 
If the leading coefficient of $\varphi$ does not depend
on $h$ this subvariety is empty. Otherwise its codimension in
$G_{d_0}$ is $l$ when $l<d_0$.
Replace $\Gamma ,f , {\cal A}(l)$ by
$ \Gamma_p ,f_p , {\cal A}_p(l)$  and consider
${\cal B}(l) = \bigcup_{p \in P} {\cal A}_p$. Then
${\cal B}(l)$ is a closed algebraic subvariety in $G_{d_0}$
whose codimension is at least $l-k$. In particular, a generic polynomial
$h^0$ from $G_{d_0}$ does not belong to ${\cal B}(k+1)$.
Let $h,\tilde {h} \in G_{d_0}$
be such that $h$
has the same $k+1$ leading coefficients as
$\tilde {h}$ does. 
Note that for every $p \in P$ and $l \leq k$ we have
$h \in {\cal A}(l)$ iff
$\tilde {h} \in {\cal A}(l)$.
Hence for every $p \in P$ and $l \leq k$
if $h^0 \in {\cal A}_p (l) -{\cal B}(k+1)$ then
$h^0 + g \in {\cal A}_p (l)- {\cal B}(k+1)$ for every $g \in G_{m}$
which is the desired conclusion. \\
\qed

Thus  suppose further that the 
function $L_{g} \cdot \Gamma $ is constant on $G_{m}$, i.e 
the points from $L_g \cap \Gamma$
do not go to infinity when $g$ runs over $G_{m}$. Note also that if $g$ 
is a small perturbation of $g^{0}$ in $G_{m}$ and 
$L_{g^{0}}$ meets $\Gamma $ normally at some point 
then $L_{g}$ meets $\Gamma $ normally at a nearby point. 
Therefore, in order to construct a subspace 
$V$ as in Lemma 4.5 we 
should take care of the set 
$\overline {w} = \{ w_{1} ,..., w_{l} \}$ at which $\Gamma$ 
meets $L_{g^{0}}$ non-normally. Let $k_i$ be the local intersection
number of $\Gamma$ and $L_{g^0}$ at $w_i$ and let
$k = k_1 + \ldots + k_l$. \vs

{\bf Definition 4.7} We call $k$ the defect of $L_{g^0}$
relative to $\Gamma$. \vs

Let $\mu_w$ be the multiplicity of $\Gamma$
at a point $w \in \Gamma$.
Put $\mu ( \Gamma )= max  (\mu_w \mid w \in \Gamma )$ and
$\alpha (\Gamma ) : = max (4 \log \mu ( \Gamma), 1)$. \vs

{\bf Proposition 4.8} {\sl Let
$\alpha : = \alpha ( \Gamma)$,
and let $S_{k }(\Gamma )$ be the subset of $G_{m}$ such that
for each $g \in S_k (\Gamma )$ the defect of $L_g$ relative to $\Gamma$ is 
at least $k$. Suppose that $k \leq m$. Then
the following are true.

(i) $S_{k }(\Gamma )$ is an algebraic subvariety of $G_{m}$
which depends algebraically on $f$.

(ii) The codimension of
$S_{k }(\Gamma )$ in $G_{m}$ is more than $( \log k)/\alpha $.

(iii)
For every $g^0 \notin S_k(\Gamma )$
there exists a nonzero polynomial $\omega$ 
in one variable such that $deg \, \omega <m$
and every generic $\C$-curve $L_g$ from $V$ described in Remark 4.4
meets $\Gamma$ at the same number of points
as $L_{g^0}$. Furthermore, if $k< m+1$ one can suppose that
$\omega (0) = 0$.

(iv) if $\Gamma$ does not contain
lines parallel to the $y$-axis then
$L_g$ is $\Gamma$ compatible for every $g \notin 
S_k (\Gamma )$.}

{\sl Proof.} Consider $L^h = \{ y+h(x)=0 \}$ for every $h \in 
G_{d_0}$. Let $w^0 \in \Gamma$
and let  $1<k^0 <d_0$.
We denote by $T_0 (w^0,k^0)$ the subset of $ G_{d_0}$
such that for every $h$ from this subset
the local intersection number of $L^h$ and $\Gamma$
at $w^0$ is at least $k^0$. Our first aim
is to estimate the codimension of $T_0(w^0,k^0)$ in
$G_{d_0}$. Let $x^0$ be
the $x$-coordinate of $w^0$ and let  $\varphi (x, h) = f(x, -h(x))$.
We denote by $\varphi_s$ the $s$-th derivative of $\varphi$
with respect to $x$. Note that $h \in {T_0}(w^0,k^0)$ iff
$$ (x^0,h(x^0))=w^0 \eqno (4.1.0)$$
$$ \varphi_s (x^0, h)=0 \eqno (4.1.s)$$
where $s=1, \ldots , k^0-1$. Let $\mu_0$ 
be the multiplicity of $\Gamma$ at $w^0$.
Then the equations (4.1.s) holds automatically for
$1\leq s \leq \mu_{0}-1$.
Let $h(x)= \sum_i c_ix^i$.
One can rewrite equations (4.1) in the form of polynomial
equations on the coefficients $\{ c_i \}$
$$\Phi_0 ( c_0 )=0, \eqno (4.2.0)$$
$$\Phi_s (\{ c_i \} ) = 0 \eqno (4.2.s)$$
where $s= \mu_{0}, \ldots , k^0-1$.
Without loss of generality we can suppose that $x^0 = 0$.
Then (4.2.0.) means that $c_0$ coincides with the
$y$-coordinate of $w^0$, and $\Phi_s$
is just the coefficient before the monomial $x^s$
in the polynomial $f(x, -h(x))$. Replacing
$y$ by $y - c_0$ we can suppose that $w^0$
is the origin. Then $c_0=0$ and one can see that
the codimension of $T_0(w^0,k^0)$ in $G_{d_0}$
is one more than the codimension of
the affine algebraic variety given by
the equations $(4.2.s)$ which do not contain
now the variable $c_0$. Consider three cases.

Case 1: $k^0=\mu_{0}\geq 2$.
The codimension of $T_0(w^0,k^0)$ in $G_{d_0}$
is $1> (\log k^0 )/\alpha $.

Case 2: $k^0 > \mu_{0} \geq 2$.
The polynomial $f$ does not
contain monomials of degree less then $\mu_{0}$.
Furthermore, $(x)^{\mu_{0}}$ cannot be the only
monomial of degree $\mu_{0}$ in $f$ with a nonzero
coefficient since $k^0 > \mu_{0}$. 
Hence $f$ contains monomials of type $y^r x^{\mu_0 -r}$ (where $r>0$)
with nonzero coefficients. Suppose also that 
$r_0$ is the maximum among such $r$'s.
Consider in $(4.2)$ the equations with
$s =  \mu_{0} t$ for some natural $t$.
The assumption on the monomials
of degree $\mu_{0}$ in $f$ implies that $\Phi_s =
\lambda (c_t)^{r_{0}} + {\Psi }_s$ where $\lambda$
is a nonzero constant, the degree of the polynomial
${\Psi }_s$ in variable $c_t$ is less than $r_{0}$.
Denote by $E_n$ the equation $(4.2.s)$ with
$s=(\mu_0)^n$ where $n$ can be any number from 1 to
$[\log (k^0-1) /\log \mu_0 ]$ ($[a]$
is the entire part of $a$). We saw already
that $E_n$ depends on $c_t$ with $t=(\mu_0)^{n-1}$. On the other hand
$E_n$ does not depend on $c_j $ where $j=(\mu_0)^i$ and $i\geq n$.
Hence the codimension of $T_0(w^0,k^0)$ in $G_{d_0}$ is
$[(\log (k^0-1) )/\log \mu_0 ] >
(\log (k^0-1))/(2 \log \mu_{0}) > 
(\log k^0 )/(4 \log \mu_0 ) \geq (\log k^0) /\alpha$.

Case 3: $\mu_{0} =1$. One can check
that for every $s$ the equation $(4.2.s)$ does not contain
$c_i$ when $i >s$ and it contains only a linear
term with $c_s$. Hence
the codimension of $T_0(w^0,k^0)$ in $G_{d_0}$
is $k^0$. That is, it is at least
$1+ (\log k^0)/\alpha $ since $k^0 \geq 2$.

Let $G$ be the subset of $G_{d_0}$ that consists of
all $h$ of form $h=h^0+g$ where $h^0$
is fixed and $g$ runs over $G_{m}$. Let
$T^0(w^0,k^0)={T_0}(w^0,k^0) \cap G$. We need to find
the codimension of $T^0(w^0,k^0)$ in $G$ in the
case when $k^0 <m$. 
Suppose that $g(x)= \sum_i b_ix^i$. 
Note that there is one-to one
correspondence between coefficients $\{ c_i \}$
and $\{ b_i \}$ where
$i=0, \ldots ,m$. Since $\Phi_s$ does not contain
unknowns $c_i$ with $i>k^0$ we see that
the codimension of $T^0(w^0,k^0)$ in $G$ is
at least $(\log k^0)/\alpha $ for singular $w^0$ and 
it is at least $1+(\log k^0)/\alpha $ when $w^0$ is regular.

Let $\w = (w_1, \ldots , w_l)$ be different points on $\Gamma$
and let $\k = (k_1, \ldots ,k_l)$ where each $k_i \geq 2$.
Consider $R_0(\w, \k ) = \bigcap^l_{i=1} T^0(w_i,k_i)$. 
Let $\hat {w} = (w_1, \ldots , w_n)$ where $n \leq l$.
Suppose that all coordinates of $\hat {w}$ are different
singular points on $\Gamma$.
Put $R^0 (\hat {w}, \k)
= \bigcup_{\w }R_0(\w ,\k )$ where
the first $n$ coordinates of $\bar {w}$ are fixed and
coincide with $\hat {w}$, and the last $l-n$ coordinates
run over all $(l-n)$-tuples of different regular points
on $\Gamma$. 
By construction,
$R_0(\w ,\k )$ depends algebraically on $\w$ and $f$.
Hence $R^0(\hat {w} ,\k )$ depends algebraically on $f$.
Since $S_k(\Gamma )$ is a union of a finite number
of sets of type $R^0(\hat {w} ,\k )$ 
(with $k=k_1 + \ldots +k_l$) this yields (i).

Let $x(i)$ be the $x$-coordinate of $w_i$.
Then the equations (4.2) for $w^0=w_i$
can be viewed as some equations on the $(k_i-1)$-jet of
$g$ at $x(i)$
and they do not impose any restrictions
on higher derivatives of $g$.
The existence of Lagrange polynomials
implies that when $k \leq m$ the codimension of
the $R_0(\w , \k )$ in $G$ is at least
$\sum_i (\log k_i)/ \alpha  + (l-n)$. Since $k_i \geq 2$ we
have $k_1 + \ldots + k_l \leq k_1 \cdots k_l$.
Hence the codimension of $R^0(\hat {w} ,\k )$ in $G$ is at least
$(\log k)/  \alpha $ which is $(ii)$.

Let $g^0 \notin S_k(\Gamma )$. Suppose that $L_{g^0}$ meets
$\Gamma $ non-normally at the points from $\w$ only.
By the Lagrange theorem, there exists
a nonzero polynomial $\omega$ of degree $k-l$
such that it has zero of order $k_i -1$ at each point $x(i)$.
Furthermore, if we allow the degree of
$\omega$ to be $k-l+1$ then we can also suppose that
that $\omega (0)=0$.
Consider $g:=g^0 +\tilde {g}\omega \in G_{m}$. 
Note that for generic $\tilde {g}$ the curve $L_g$
meets $\Gamma $ non-normally at the points of $\w$ only.
In particular, it meets $\Gamma$ at the same number of points
as $L_{g^0}$ does which is (iii). 
(Recall that due to Lemma 4.6
we suppose that the intersection numbers $L_g \cdot \Gamma$ and
$L_{g^0} \cdot \Gamma$ are the same.)
Now (iv) follows from Lemma 4.5. \\
\qed

{\bf Corollary 4.9} {\sl Let
$\hat {f}$ be a polynomial on $\C^{3}$. For every
$g \in G_{m}$ consider the curve $\Gamma (g)$
that is the zero locus of the polynomial
$f^0 (x,y ) = \hat {f} (x,y,g(x))$. Let $f(x,y ) =
\hat {f} (x , y, -y-h^0(x))$ where
$h^0$ is as before Lemma 4.3.
Suppose that $\Gamma$ is the zero locus of $f$
and $S_k(\Gamma )$ is the same as in Proposition 4.8.
If $k <m$ then

(i)
for every $g^0 \notin S_k(\Gamma )$
there exists a nonzero polynomial $\omega$ 
in one variable such that $deg \, \omega <m$
and every generic $\C$-curve $L_g$ from $V$ described in Remark 4.4
meets $\Gamma ({g^0})$ at the same number of points
as $L_{g^0}$;

(ii) if $\Gamma ({g^0})$ does not contain
lines parallel to the $y$-axis then
$L_g$ is $\Gamma ({g^0})$ compatible for every $g \notin 
S_k (\Gamma )$.}

{\sl Proof.}  Put $\psi (x) = f(x, -h(x))$ 
where $h(x) = h^0(x) + g^0(x)$. The roots
of $\psi$ are the $x$-coordinates of the
points where $L_{g^0}$ meets $\Gamma$. The points
where $L_{g^0}$ meets $\Gamma$ non-normally correspond to
the multiple roots of $\psi$ and the local
intersection number of $L_{g^0}$ and $\Gamma$ at each
of these points  coincides with the corresponding multiplicity.
Note that
$f^0 (x,y ) = f (x,y )+ \tilde {f} (x,y )$
where $\tilde {f}$ belongs to the ideal
generated by the polynomial $y + h(x)$.
Since $y+h(x) \equiv 0$ on $L_g$ we have $\psi (x) =f^0(x, -h(x))$.
This implies that $L_{g^0}$ meets $\Gamma ({g^0})$ at the
same points as $L_{g^0}$ meets $\Gamma$ and with the
same multiplicity. In particular, the defects
of $L_{g^0}$ relative to $\Gamma$ and relative to $\Gamma ({g^0})$
coincide. Let $x(i)$ be the $i$-th multiple
root of $\psi$ where $i=1, \ldots , l$ and let
$k_i$ be the multiplicity of this root.
Suppose that $w_i \in \C^{2}$ is the point
with $x$-coordinate $x(i)$ where $L_{g^0}$ meets $\Gamma ({g^0})$.
As we showed in the proof of Proposition 4.8
$L_g$ meets $\Gamma ({g^0})$ at $w_i$
with local intersection number at least $k_i$
when $g$ has the same $(k_i-1)$-jet at $x(i)$ as
$g^0$ has. Since $k=k_1+ \ldots +k_l <m$ there exists
a nonzero polynomial $\omega (x)$ of degree less than $m$ which has
zeros at $x(i)$ of multiplicity $k_i$ for every $i$.
Consider generic $g(x)= g^0 (x)+
\omega (x) \tilde {g} (x)=0$ where 
$\tilde {g} \in G_{n}$ and 
$ n=m-deg \omega$. By Lemma 4.6,
the intersection number $L_g \cdot \Gamma ({g^0})$ does not
depend on $g$, and since $g$ is generic
the number of points where $L_g$ meets $\Gamma ({g^0})$
normally is at least the same as the number
of points where $L_{g^0}$ meets $\Gamma ({g^0})$ normally.
This implies (i). Lemma 4.5 yields (ii). \\
\qed

{\bf Corollary 4.10} {\sl Let 
$\{ \hat {f}_p \mid \, p \in P \}$ be a family of polynomials
on $\C^{3}$. For every
$g \in G_{m}$ consider the curve $\Gamma_p(g)$
that is the zero locus of the polynomial
$f^0_p (x,y ) = \hat {f}_p (x,y,g(x))$. Let $f_p(x,y ) =
\hat {f}_p (x,y,-y -h^0(x))$ where
$h^0$ is as before Lemma 4.3.
Let $\Gamma_p$ be the zero locus of $f_p$. Suppose that
$\theta = \max_p \alpha ( \Gamma_p )$ where $\alpha (\Gamma_p)$ is
as before Proposition 4.8.
Suppose that 
$m> \exp (\theta  dim \, P )$. Then 

(i) the closure
$S(P) $ of ${\ds \bigcup _{p\in P}} S_m (\Gamma_{p} )$ in 
$G_{m}$ is a proper algebraic subvariety
where $S_m (\Gamma_p)$ has the same meaning as in Proposition 4.8;

(ii) for every $g^0 \notin S(P)$ and every $p \in P$
there exists a nonzero polynomial $\omega_p$ 
in one variable such that for
a generic element $g \in G_{m}$ of
form $g= g^0 + \tilde {g} \omega_p$ the curve
$L_g$ meets $\Gamma_p(g^0)$ at the same number of points as $L_{g^0}$ does;

(iii) if $g^0 \notin S(P)$ and none of the curves
$\{ \Gamma_p (g^0) \}$ contains 
lines parallel to the $y$-axis
then $L_{g^0}$
is compatible relative to the family $\{ \Gamma_{p}(g^0) \}$.} \vs

{\em Proof.}
The codimension of $S(P)$ in
$G_{m}$ is more than $(\log m ) /\theta  -dim \, P >0$, by Proposition 4.8. 
Hence $G_{m} -S(P)$ is Zariski open which is the first statement
of this Corollary. The second and the third
statements follow from Corallary 4.9. \\
\qed

{\bf 5. Main result. } \vs

We shall consider first the three-dimensional case. \vs

{\bf Proposition 5.1.} {\em Let $\{ H_{p} \}$ be a family of hypersurfaces in 
$\C^{3}$ with parameter $p\in P$. Then 
there exists a coordinate system in $\C^3$
such that some plane in this system
is strictly compatible relative to $ \{ H_{p} \}$. }

{\sl Proof.}
Let $(x, y, z)$ be a coordinate system in $\C^3$ such that
none of the surfaces $H_p$ contains a plane $z = const$
(it is enough to require that the restriction of the
projection $(x,y,z) \to (x,z)$ to every $H_p$ is
finite which can be done by Lemma 4.1).
Let $\varphi$ be a polynomial in one variables such that
$deg \, \varphi
>>  \max_p (deg \, f_p, dim \, P)$. Put $q=(p, \varphi )$
and $f^q(x,y,z)=f_p(x+\varphi (z),y,z)$
where $f_p$ is a defining polynomial for $H_p$.
Suppose that $Q$ is the variety of all $q$'s and
$H^q$ is the zero locus  of $f^q$.
Then $\{ H^q \mid \, q \in Q \}$ is a family of
hypersurfaces which is invariant under automorphisms of form
$(x,y,z) \to (x + \varphi (z), y,z)$. 
If $Q_{0}$ is a subvariety of $Q$ that consists of
$q=(p,\varphi_0 )$ for some fixed $\varphi_0$
then the subfamily $\{ H^q \mid \, q \in Q_{0} \}$ coincides
with $\{ H_p \mid \, p \in P \}$ after a polynomial
coordinate substitution.

Let $h$ be a polynomial in one variable such that
$deg \, h >> deg \, \varphi$. Replace
$(x,y,z)$ by $(x, y+h(x), z)$. Then the equation of $H^q$
becomes $f^q(x,y+h(x),z)=
f_p(x + \varphi (z), y+h(x+\varphi (z)),z)=0$.
By Lemma 4.1 
the restriction of the projection
$\rho : \C^3 \to \C^2$ (given by $(x,y,z) \to (x,y)$)
to every $H^q$ is finite.
When $h$ is fixed
the image of the ramification set of this restriction
is a curve $\Gamma^q$ in the $(x,y)$-plane which depends on
$q \in Q$. If we replace $h$ by $h+g$
then each curve $\Gamma^q$ in the $(x,y)$-plane
must be replaced by its image under automorphism
$(x,y) \to (x, y+g(x))$. Suppose that $deg \, g << deg \, h$
but still $deg \, g >> dim \, Q$. By Lemma 4.6 one can suppose
that for every $q$ the intersection
number of $\Gamma^q$ and the curve $L_{g}$ given by
$y+g(x)=0$ is finite and does not depend on $g$. In particular,
$L_{g}$ is not a component of $\Gamma^q$. 
Hence condition (2) from Definition 3.1 is true for $(x,z)$-plane
$A_0$ with respect to the family of surfaces $\{ H^q | \, q \in Q \}$.

Consider the family of
curves $\Theta_q$ in $A_0$ which are
the intersections $H^q \cap A_0$. That is,
the equation of $\Theta_q$ is
$f_p(x + \varphi (z), h(x+ \varphi (z)),z)=0$.
Show that none of these curves contain a line
$z=c$. Consider the curves $\{ \Lambda_{c,q} \}$
in the $(x,y)$-plane given by the equations
$f_p(x + \varphi (c), y+h(x+ \varphi (c)) ,c)=0$.
Note that $\Theta_q$ contains the line $z=c$
iff $\Lambda_{c,q}$ contains the $y$-axis.
But $\Lambda_{c,q}$ cannot contain the $y$-axis
since the restriction of the projection
$(x,y) \to y$ to every $\Lambda_{c,q}$ is finite
by Lemma 4.1.

If $\Sigma$ is a curve in $\C^2$
denote by $\mu (\Sigma )$ the maximum of 
multiplicities of its points.
Consider the curve $\Sigma$ in $A_0$ given by the equation
$f_p(x, h(x),z)=0$. Note
that $\mu (\Sigma ) $ is at most
$deg_z \, f_p(x, y,z)
\leq deg \, f_p$
since this is the number of points at which
a generic line $x= const$ meets $\Sigma $.
Since $\Theta_q$ can be obtained from $\Sigma$
by an automorphism we see that
$\max_{q \in Q} \mu (\Theta_q)$ is bounded by
$\max_p deg \, f_p << deg \, \varphi$.
Hence $\max_{q \in Q} \alpha (\Theta_q)
<< deg \, \varphi$ where $\alpha$ is the same
as before Proposition 4.8.
Corollary 4.10 implies that
one can choose $\varphi =\varphi_0$ so that the $z$-axis in $A_0$
is $\Theta_q$ compatible for every $q$ of form
$(p, \varphi_0)$ where $p$ is arbitrary.
Hence condition (1) from Definition 3.1 is true for
$A_0$ with respect to the subfamily $\{ H^q \mid q \in Q_0 \}$.

Let $\Gamma$ be the
image in the $(x,y)$-plane of the ramification set
of the restriction of $\rho$ to the surface given by
$f^q(x,y,z)=0$.
Note that $\Gamma^q$ 
can be obtained from $\Gamma$
by an automorphism of the $(x,y)$-plane and, therefore,
$\alpha (\Gamma_q ) = \alpha (\Gamma )$. Hence
$M= max_{q \in Q} \alpha (\Gamma^q )$ is independent from
$h$. Since we have no restrictions on the degree of $h$
we can suppose from the beginning that $deg \, h$
is much greater than $\exp (M dim \, Q)$.
Let $deg \, g < deg \, h$ but still
$deg \, g >> max (deg \, \varphi_0 , \exp (M, dim \, Q))$.
By virtue of Proposition 4.8 one can suppose that
the following are true:
for every $q \in Q$ there exists a nonzero polynomial $g$ such that
$g(0)=0$ and the line $L_0$ (i.e., the $x$-axis)
meets $\Gamma^q$ at the same number of points 
(counting without multiplicity) as
the curve $L_{cg}$ for generic $c \in \C$.
This gives conditions (3) from Definition 3.1
for $A_0$ with respect
to the family $\{ H^q \mid \, q \in Q \}$.
That is, $A_0$ is strictly compatible
relative to the family $\{ H^q \mid \, q \in Q_0 \}$
by Theorem 3.2.\\
\qed

{\bf Remark 5.2} Using the fact that
$deg_z \, f_p(x + \varphi (z), h(x+ \varphi (z)),z)=
deg_z \,  f_p(x + \varphi (z),y+ h(x+ \varphi (z)),z)$,
it is not difficult to check that every plane
from a Zariski open neighborhood of $A_0$ in the variety
of planes (i.e., every generic plane)
is strictly compatible
with respect to the subfamily $\{ H^q \mid q \in Q_0 \}$.
In fact, one can construct more sophisticated
coordinate substitutions such that every plane becomes
strictly $H_p$ compatible for every $p \in P$.
The same remark is applicable in the
case of an arbitrary dimension. \vs

{\bf Lemma 5.3} {\sl Let
$H$ be a hypersurface in $\C^{n+1}$, 
let $\lambda$ be a coordinate function on $\C^{n+1}$,
and let $R$ be a hyperplane in $\C^{n+1}$ which is not given
by $\lambda = const$.
For $c\in \C$ put $C_c =\lambda^{-1} (c), \;\; 
R_c =C_c \cap R, \;\; H_c =H\cap C_{c}$. Suppose that for a 
generic $c\in \C$ the manifold $R_{c}$ is strictly $H_{c}$ 
compatible in $C_{c} \cong \C^{n}$ and that for every $c \in \C$
the set $R_c$ is not an irreducible component of $H_c$. 
Then $R$ is 
strictly $H$ compatible. }\vs

{\sl Proof.} Let $S$ be a finite set in $\C$. 
Put $R(S)=R-(\lambda^{-1} (S)\cup H), \; 
C(S)=\C^{n+1} -(\lambda^{-1} (S) \cup H), 
\;\; \lambda_{S} =
\lambda \vline _{C(S)}$, and $\lambda'_{S} =\lambda \vline _{R(S)}$.
Choose the finite set $S \subset \C$ so that the mappings 
$\lambda _{S} :C(S)\rightarrow \C -S$ and 
$\lambda'_{S} :R(S)\rightarrow \C -S$ 
are fibrations. Let $F$ be the fiber of 
$\lambda_{S}$ and let $F'$ be the 
fiber $\lambda'_{S}$. Then, by assumption, 
the natural embedding $e:F'\rightarrow F$ generates an isomorphism 
$e_{*} :\pi _{1} (F')\rightarrow \pi _{1} (F)$. 
We have also two other embeddings 
$i:F\rightarrow C(S)$ and $i':F'\rightarrow R(S)$ 
which implies the commutative diagram 

\[ \begin{array}{lcccr} 0 \rightarrow &
\pi _{1} (F) & \stackrel{i_{*}}{\rightarrow } \pi _{1} (C(S)) \stackrel{\lambda_{S_{*}}}{\rightarrow } &
\pi _{1} (\C-S) & \rightarrow 0 \\
\uparrow & \uparrow e_{*} & \uparrow & \uparrow id & \uparrow \\
0 \rightarrow & \pi _{1} (F') & \stackrel{i'_{*}}{\rightarrow } \pi _{1} (R(S)) \stackrel{\lambda'_{S_{*}}}{\rightarrow } &
\pi _{1} (\C-S) & \rightarrow 0. \end{array} \]

\noin The five isomorphisms lemma implies that $\pi _{1} (C(S))$ and $\pi _{1} (R(S))$ are isomorphic. 
Since $R_s$ is not contained in $H_s$ for every $s$ one can
choose simple loops $\gamma _{s}$ in $R(S)$ around 
each hypersurface $R_s$ with $s\in S$ such that 
$\gamma _{s}$ is contractible in 
$R-H$. Consider the natural embedding $j':R(S)\rightarrow R-H$. It 
generates an epimorphism $j'_{*} :\pi _{1} (R(S)) 
\rightarrow \pi _{1} (R-H)$ and, obviously, 
$[\gamma _{s} ]\in ker \, j'_{*}$ for every $s\in S$. Moreover, if $N'$ is the smallest normal subgroup 
in $\pi _{1} (R(S))$ that contains all $[\gamma _{s} ], \; s\in S$, then $ker 
\, j'_{*} =N'$, which follows from two simple geometrical observations: 
\vs

-each two-cell in $R-H$ becomes transversal to $\lambda^{-1} (S)\cap R$ after a perturbation, i.e. every 
contractible loop in $R-H$ can be viewed as a product of simple 
contractible loops around hyperplanes $R_s, \; 
s\in S$; \vs

-every simple contractible loop of this type is conjugate to some $[\gamma _{s} ]$ as an element of $\pi _{1} (R(S))$. \vs

Similarly, we can consider the embedding 
$j:C(S) \rightarrow \C^{n} -H$. It generates an epimorphism $j_{*} :\pi _{1} 
(C(S))\rightarrow \pi _{1} (\C^{n} -H)$. Then $ker \, j_{*}$ coincides with the smallest normal subgroup $N$ of $\pi _{1} 
(C(S))$ that contains all $[\gamma _{s} ], \; s\in 
S$, where we treat $[\gamma _{s} ]$ as an element of $\pi _{1} (C(S))$ now. 
This yields an isomorphism between $\pi _{1} (R-H)$ and $\pi _{1} (\C^{n} -H)$ which concludes the proof. \\
\qed

{\bf Theorem 5.4} {\em Let $\{ H_p \}$ be a family of 
hypersurfaces in $\C^{n}$ with parameter $p\in P$. Then there
exists a polynomial coordinate system in $\C^{n+1}$ such that
some plane is strictly $H_{p}$ 
compatible for every $p\in P$.}\vs

{\sl Proof.} We shall use induction. The first step of 
induction is Proposition 5.1. Assume that for every family of hypersurfaces
in $\C^{n}$ there exist a coordinate system such that
some hyperplane in this system 
is strictly compatible relative to this family.
Consider a family of hypersurfaces $\{ H_{p} \}$ in $\C^{n+1}$.
Let $\x = (x_1, \ldots , x_{n+1})$ 
be a coordinate system in $\C^{n+1}$ and let
$\lambda$ coincide with $x_1$ on $\C^{n+1}$. We can choose 
$\x$ so that none of $H_p$ contains a hyperplane $x_1 =const$
since the restriction of the projection $\x \to (x_1, \ldots ,
x_n)$ to $H_p$ can be supposed to be finite, 
by the analogue of Lemma 4.1 in the case of higher dimensions.
Put $C_c = \lambda^{-1} (c)$.
We can view $H_{p,c} =H_{p} \cap C_{c}$ as a hypersurface in 
the fiber $C_c$.
Put $Q=P\times \C$. Then we can consider 
$\{H_{p,c}= H^{q} \}$ with $q=(p,c)\in Q$ as 
a family of hypersurfaces in $\C^{n}$ . By induction, a coordinate 
system $\bar {y}= (y_2, \ldots , y_{n+1})$ 
in $\C^{n}$ can be chosed so that some
hyperplane $E$ in $\C^n$ is strictly compatible relative
to $\{ H^q \}$. In particular, none of
the hypersurfaces $H^q$ contain $E$.
Let $R$ be the hyperplane $\tau^{-1} (E)$
in the coordinate system $(x_1, \bar {y} )$ in $\C^{n+1}$
where $\tau$ is the natural projection $(x_1, \bar {y}) \to \bar {y}$.
By Lemma 5.3, the hyperplane $R$ is strictly compatible relative 
$ \{ H_p \}$.
Therefore, we can reduce 
dimension by induction, which implies our Theorem. \\
\qed
\vspace{.13in}

\begin{center} {\bf Bibliography.} \end{center} 

\noin [1] O. Zariski, {\em A theorem on the Poincar$\acute{e}$ group of 
an algebraic hypersurface}, Annal of Math., {\bf 38}(1937), 131-141. \vs

\noin [2] Sh. Kaliman, {\em On the Jacobian conjecture}, Proc. AMS, 
{\bf 117}(1993), 45-51. \vs

\noin [3] E. van Kampen, {\em On the connection between the fundamental
groups of some related spaces}, Amer. J. of Math., {\bf 55}(1933),
261--267. \vs

\end{document}